\documentclass[11pt]{article}
\usepackage{mathrsfs}
\usepackage{amssymb}
\usepackage{amsmath,amsthm,amssymb,mathrsfs,amsbsy,bm}
\usepackage{graphicx}
\usepackage{epsfig}
\usepackage{enumerate}
\usepackage{color, xcolor} 

\usepackage[colorlinks=true, allcolors=blue]{hyperref}

\usepackage[letterpaper,top=2cm,bottom=2cm,left=3cm,right=3cm,marginparwidth=1.75cm]{geometry}

\numberwithin{equation}{section}

\newtheorem{theorem}{\bf{Theorem}}[section]
\newtheorem{lemma}{\bf {Lemma}}[section]
\newtheorem{define}{\bf{Definition}}[section]
\newtheorem{corollary}{\bf{Corollary}}[section]
\newtheorem{proposition}{\bf{Proposition}}[section]
\newtheorem{remark}{\bf{Remark}}[section]

\newcommand{\mbE}{\widehat{\mathbb{E}}}
\newcommand{\mbe}{\widehat{\mathcal{E}}}
\newcommand{\V}{\mathbb{V}}
\newcommand{\mv}{\mathcal{V}}

\newcommand*{\dif}{\mathop{}\!\mathrm{d}}

\title{Marcinkiewicz-type laws of large numbers for pseudo-independent random variables under sublinear expectations}
\author{Jialiang Fu\\ \sl \small Academy of Mathematics and Systems Science, Chinese Academy of Sciences, \\ \sl \small Beijing, China\\ \sl \small E-mail: fujialiang@amss.ac.cn }
\date{}

\begin{document}

\maketitle

\begin{abstract}
\par As a kind of independence of random variables under sublinear expectations, pseudo-independence is weaker than Peng's independence. We shall give Marcinkiewicz-type weak and strong laws of large numbers for pseudo-independent random variables under the framework of sublinear expectations.   
\end{abstract}

\section{Introduction}
 \ \ \ \ \  Let $\{X_n, n\geq 1\}$ be a sequence of independent and identically distributed (i.i.d.) random variables in $(\Omega,P,\mathcal{F})$ and $S_n\triangleq\sum_{i=1}^{n}X_i$. The strong law of large numbers (SLLN) established by Kolmogorov states that if $E_{P}[|X_1|]<\infty$, then
 \begin{equation}
 	P\left(\lim_{n\rightarrow\infty}\frac{S_n}{n}=E_{P}[X_1]\right)=1,\notag
 \end{equation}
 where $E_{P}$ is the expectation with respect to the probability measure $P$. As the generalization of Kolmogorov's SLLN, the classical Marcinkiewicz-type SLLN states that if $E_{P}[|X_1|^{r}]<\infty$, for some $r\in[1,2)$, then
 \begin{equation}
 	P\left(\lim_{n\rightarrow\infty}\frac{S_n-nE_{P}[X_1]}{n^{\frac{1}{r}}}=0\right)=1.\notag
 \end{equation}
 It's clear that Marcinkiewicz-type SLLN gives the $O\left(\frac{1}{n^{1-1/r}}\right)$ convergence rate of Kolmogorov's SLLN when the random variables have $r$-order moments for some $r\in(1,2)$.
 \par The nonlinear expectation is a powerful tool when the probability models have uncertainty. A common nonlinear expectation is the sublinear expectation, which can be represented as the supremum of a set of linear expectations. Corresponding to classic linear expectations, the law of large numbers (LLN) can also be established under the sublinear expectation. Peng (2007) \cite{Pen19} introduced the concepts of independence and identical distribution for random variables under sublinear expectations, and established the weak law of large numbers (WLLN) for i.i.d. random variables under sublinear expectations. Extensive literature has been developed about the law of large numbers under sublinear expectations. Chen \cite{Ch16} established a strong form of the LLN for i.i.d. sequences. Song \cite{S23} obtained a strong law of large numbers under sublinear expectations and studied the triviality of the tail $\sigma$-algebra of a sequence of i.i.d. random variables under a sublinear expectation. Zhang \cite{Zha24} studied the limit points of the strong law of large numbers under sublinear expectations. It's necessary to extend the LLN for various kinds of dependent random variables because real-world data can be dependent and more complex. Under the framework of sublinear expectations, Zhang \cite{Zha16,GZ24} obtained the strong laws of large numbers for negatively dependent and identically distributed random variables and for $m$-dependent and stationary random variables. Fu \cite{Fu25} obtained the LLN for blockwise $m$-dependent random variables under sublinear expectations. The concept of pseudo-independence under sublinear expectations was introduced by Guo and Li \cite{GL21}, and they established the weak and strong laws of large numbers for pseudo-independent random variables with finite 1-order Choquet integrals. In this paper, we shall give Marcinkiewicz-type weak and strong laws of numbers for pseudo-independent random variables with finite $r$-order Choquet integrals ($r\in [1,2)$).
 
 \par The structure of this paper is as follows. Section 2 will present some preliminaries for the sublinear expectations. Section 3 will show Marcinkiewicz-type weak and strong laws of numbers, and their detailed proofs are given in Section 4.

\section{Preliminaries}
\ \ \ \ \ In this section, we shall present some basic concepts of sublinear expectation theory and some existing results. One can refer to \cite{Pen19} for more details on the former.
\par Let $(\Omega,\mathcal{F})$ be a given measurable space and $\mathcal{M}$ be the collection of all probability measures on $(\Omega,\mathcal{F})$. $X\in \mathcal{F}$ denotes that $X$ is a random variable defined on $(\Omega,\mathcal{F})$. For a given subset $\mathcal{P}\subseteq\mathcal{M}$, the upper expectation and the lower expectation with respect to $\mathcal{P}$ are defined as follows:\\
\begin{equation}
	\mbE[X]\triangleq\sup_{P\in \mathcal{P}}E_{P}[X],\ \mbe[X]\triangleq-\mbE[-X]=\inf_{P\in\mathcal{P}}E_{P}[X],\notag
\end{equation}
for $X \in \mathcal{F}$ with $\mbE[X]$ and $\mbe[X]$ being finite.
\par $\mbE$ is a sublinear expectation which satisfies:
\begin{itemize}
\item[(i)] Monotonicity: $\mbE[X]\leq\mbE[Y]$ if $X\leq Y$;
\item[(ii)] Constant preserving: $\mbE[c]=c$ for $c\in \mathbb{R}$;
\item[(iii)] Sub-additivity: $\mbE[X+Y]\leq\mbE[X]+\mbE[Y]$;
\item[(iv)] Positive homogeneity: $\mbE[\lambda X]=\lambda\mbE[X]$ for $\lambda\geq0$.
\end{itemize}
The upper probability and the lower probability are defined as follows:
\begin{equation}
	\V(A)\triangleq\sup_{P\in\mathcal{P}}P(A),\quad \mv(A)\triangleq\inf_{P\in\mathcal{P}}P(A),\quad \forall A \in \mathcal{F}.\notag
\end{equation}
It's clear that $\V$ and $\mv$ are conjugate to each other, which means 
\begin{equation}
\V(A)+\mv(A^c)=1,\ \forall A \in \mathcal{F}.\label{eq5}
\end{equation}
\par For $X\in \mathcal{F}$, we define the Choquet integrals $(C_{\V},C_{\mv})$ by
\begin{equation}
	C_V(X)\triangleq\int_0^{\infty}V(X\geq t)\dif t+\int_{-\infty}^0[V(X\geq t)-1]\dif t\notag
\end{equation}with $V$ being replaced by $\V$ and $\mv$ respectively. It's clear that $\mbE[|X|]\leq C_{\V}(|X|)$. 
\begin{define}
Let $\mbE_1$ and $\mbE_2$ be two sublinear expectations on $(\Omega_1,\mathcal{F}_1)$ and $(\Omega_2,\mathcal{F}_2)$ respectively. A random variable $X_1$ in $(\Omega_1,\mathcal{F}_1)$ under $\mbE_1$ is said to be identically distributed with another random variable $X_2$ in $(\Omega_2,\mathcal{F}_2)$ under $\mbE_2$, denoted by $X_1\overset{d}{=}X_2$, if
$$
	\mbE_1[\varphi(X_1)]=\mbE_2[\varphi(X_2)],\quad\forall\varphi\in C_{b,Lip}(\mathbb{R}),
$$ where $C_{b, Lip}(\mathbb{R})$ denotes the set of all bounded Lipschitz functions on $\mathbb{R}$.
 A sequence $\{X_n, n\geq 1\}$ of random variables is said to be identically distributed if $X_i\overset{d}{=}X_1$ for each $i\geq 1$.

\end{define}
\begin{define}
	(Peng's independence) Let $X$ and $Y$ be two random variables in $(\Omega,\mathcal{F})$. $Y$ is said to be independent of $X$ if for each test function $\varphi\in C_{b,Lip}(\mathbb{R})$ we have $$\mbE[\varphi(X,Y)]=\mbE[\mbE[\varphi(x,Y)]\vert_{x=X}].$$
\end{define}
 Generally, it is important to observe that under the framework of sublinear expectation, $Y$ is independent of $X$ does not in general imply that $X$ is independent of $Y$, which is different from the classical linear expectation. One can check the Example 1.3.15 in \cite{Pen19} for details. A sequence of random variables $\{X_n, n\geq 1\}$ is said to be independent if $X_{i+1}$ is independent of $(X_1,\cdots,X_i)$ for each $i\geq 1$. It is easy to check that if $\{X_1,\cdots,X_n\}$ are independent, then $\mbE[\sum_{i=1}^nX_i]=\sum_{i=1}^n\mbE[X_i]$. 
 \par Let $X$ and $Y$ be two random variables in $(\Omega,\mathcal{F})$. If $X\overset{d}{=}Y$, then 
 \begin{equation}
 	\V(X\geq x+\epsilon)\leq\V(Y\geq x)\leq\V(X\geq x-\epsilon), \forall x\in \mathbb{R}, \forall \epsilon>0,\label{eq3}
 \end{equation}
 \begin{equation}
 	C_{\V}(X)=C_{\V}(Y).\label{eq4}
 \end{equation}
 Actually, we can find a $\varphi\in C_{b,Lip}(\mathbb{R})$ such that $I_{[x+\epsilon,\infty)}(y)\leq\varphi(y)\leq I_{[x,\infty)}(y)$, then 
 \begin{equation}
 	\V(X\geq x+\epsilon)=\sup_{P\in \mathcal{P}}E_{P}[I_{[x+\epsilon,\infty)}(X)]\leq \sup_{P\in \mathcal{P}}E_{P}[\varphi(X)]=\mbE[\varphi(X)],\notag
 \end{equation}
 \begin{equation}
 	\mbE[\varphi(Y)]=\sup_{P\in \mathcal{P}}E_{P}[\varphi(Y)]\leq \sup_{P\in \mathcal{P}}E_{P}[I_{[x,\infty)}(Y)]=\V(Y\geq x).\notag
 \end{equation}
 We get the first inequality in \eqref{eq3} by $X\overset{d}{=}Y$, and the second inequality is similar. By \eqref{eq3}, it's true that $\V(X\geq x)=\V(Y\geq x)$ if $x$ is a continuous point of both functions $\V(X\geq y)$ and $\V(Y\geq y)$. But we know that both of them are non-increasing functions, and a monotone function
 has at most a countable number of discontinuous points in real analysis. So
 \begin{equation}
 	\V(X\geq x)=\V(Y\geq x) \quad \text{for all but except countable many}\ x,\notag
 \end{equation} 
 and hence 
 \begin{equation}
 	C_{\V}(X)=C_{\V}(Y).\notag
 \end{equation}

 \par Guo and Li \cite{GL21} introduced the definition of pseudo-independence under sublinear expectations. The motivation for introducing pseudo-independence is inspired by the property of independence under linear expectations, which is that given a random variable $X$ on $(\Omega,\mathcal{F},P)$ and a sub-$\sigma$-algebra $\mathcal{G}\subset\mathcal{F}$,
 \begin{equation}
 	E_P[\varphi(X)|\mathcal{G}]=E_P[\varphi(X)], \ \forall\varphi\in C_{b,Lip}(\mathbb{R}) \iff X\ \text{is independent of}\ \mathcal{G}.\label{eq9}
 \end{equation} 

  \begin{define}
 	Let $\{X_n, n\geq1\}$ be a sequence of random variables on $(\Omega, \mathcal{F}, \mbE)$. Define
 	$\mathcal{F}_n\triangleq\sigma(X_1,\dots,X_n)$ and $\mathcal{F}_0\triangleq\{\emptyset, \Omega\}$. If for each $P\in \mathcal{P}$, we have
 	\begin{equation}
 		E_{P}[\varphi(X_n)|\mathcal{F}_{n-1}]\leq\mbE[\varphi(X_n)],\ P-a.s.,\ \forall \varphi \in C_{b,Lip}(\mathbb{R}),\label{eq8}
 	\end{equation}
 	then we call that $X_n$ is pseudo-independent of $\mathcal{F}_{n-1}$. For the case $n\geq 2$, we can also call that $X_n$ is pseudo-independent of $(X_1,\dots,X_{n-1})$. A sequence of random variables $\{X_n,n\geq1\}$ is called pseudo-independent if \eqref{eq8} holds for each $n\geq 1$.
 \end{define}
 \begin{remark}
\eqref{eq8} is equivalent to 
\begin{equation}
	\mbe[\varphi(X_n)]\leq E_{P}[\varphi(X_n)|\mathcal{F}_{n-1}]\leq\mbE[\varphi (X_n)],\ \forall\varphi\in C_{b,Lip}(\mathbb{R}),\ \forall P\in \mathcal{P},\notag
\end{equation} 	
which means the conditional expectation lies in a non-random interval. If $\mathcal{P}=\{P\}$ is a singleton set, then the pseudo-independence is the independence in the classical probability theory by \eqref{eq9}.
 \end{remark}
 The next is Proposition 2.4 in \cite{GL21}, which shows that Peng's independence implies the pseudo-independence.
 \begin{proposition}
 	An independent sequence of random variables $\{X_n,n\geq1\}$ on $(\Omega,\mathcal{F},\mbE)$ is pseudo-independent.
 \end{proposition}
 \begin{remark}
 	Example 3.6 in \cite{Li21} shows that pseudo-independence can not imply Peng's independence, so it's clear that pseudo-independence is weaker than Peng's independence.
 \end{remark}
 The next results are Guo and Li's in \cite{GL21}.
 \begin{theorem}
 Let $\{X_n, n\geq1 \}$ be a pseudo-independent sequence under $\mbE$ satisfying: \\
 there exist a random variable $X$ with $C_{\V}(|X|)<\infty$ and a constant $C$ such that
 \begin{equation}
 	\V(|X_n|\geq x)\leq C\V(|X|\geq x), \quad \forall x\geq0, \ \forall n\geq 1.\notag
 \end{equation} 
 Let
 \begin{equation}
 	\overline{\mu}\triangleq\limsup_{n\rightarrow\infty}\frac{1}{n}\sum_{i=1}^{n}\mbE[X_i],\ \underline{\mu}\triangleq\liminf_{n\rightarrow\infty}\frac{1}{n}\sum_{i=1}^{n}\mbe[X_i],\ S_n\triangleq\sum_{i=1}^{n}X_i.\notag
 \end{equation} 
Then
 \begin{align}
 	(WLLN)&\quad \quad  \lim_{n\rightarrow\infty}\mv\left(\underline{\mu}-\epsilon<\frac{S_n}{n}<\overline{\mu}+\epsilon\right)=1, \quad \forall \epsilon>0,\label{eq1}\\
 	(SLLN)&\quad \quad  \mv\left(\underline{\mu}\leq\liminf_{n\rightarrow\infty}\frac{S_n}{n}\leq\limsup_{n\rightarrow\infty}\frac{S_n}{n}\leq\overline{\mu}\right)=1.\label{eq2}
 \end{align}
 \end{theorem}
 \eqref{eq1} and \eqref{eq2} are Kolmogorov-type laws of large numbers, and we shall give the corresponding Marcinkiewicz-type laws of large numbers in Section 3. The next is the Borel-Cantelli lemma in \cite{Zha16}. It holds because $\V$ is a countably sub-additive capacity. 
 \begin{lemma}
 	Let $\{A_n,n\geq 1\}$ be in $\mathcal{F}$. If $\sum_{n=1}^{\infty}\V(A_n)<\infty$, then
 	\begin{equation}
 		\V(A_n,i.o.)=0,\notag
 	\end{equation} 
 	where $\{A_n,i.o.\}\triangleq\bigcap_{n=1}^\infty\bigcup_{i=n}^\infty A_i$.
 \end{lemma}
 The following lemma is a martingale convergence result in \cite{HH80}, and it plays an important role in the proof of \eqref{eq6}. 
 \begin{lemma}\label{le1}
 	Let $\{S_n\triangleq\sum_{i=1}^{n}X_i, \mathcal{F}_n, n\geq 1\}$ be a martingale and let $\mathcal{F}_0$ be the trivial $\sigma$-algebra. Then $S_n$ converges $a.s.$ on the set $\{\sum_{i=1}^{\infty}E[X_i^2|\mathcal{F}_{i-1}]<\infty\}.$
 \end{lemma}
 The next is Kronecker's lemma, which is a common tool in proving limit theorems of probability theory.
 \begin{lemma}
 	Let $\{x_n,n\in \mathbb{N}^*\}$ and $\{b_n,n\in \mathbb{N}^*\}$ be infinite sequences of real numbers, $0<b_n\uparrow\infty$. If $\sum_{n=1}^{\infty}\frac{x_n}{b_n}<\infty$, then $\lim_{n\rightarrow\infty}\frac{1}{b_n}\sum_{k=1}^{n}x_k=0$.
 \end{lemma} 
 \par  Throughout the whole paper, we denote $x\vee y\triangleq\max\{x,y\}, x\wedge y\triangleq\min\{x,y\}, x^+\triangleq x\vee0, x^-\triangleq(-x)\vee0$ for real numbers $x$ and $y$. $\mathbb{R}$ represents all real numbers. $\mathbb{N}$ represents all natural numbers, and $\mathbb{N}^*$ represents all non-zero natural numbers. $0<C$ is a constant that may change from line to line.

\section{Main results}
\ \ \ \ \ Our results are the following theorems.
\begin{theorem}\label{th1}
	 Let $\{X_n, n\geq1 \}$ be a pseudo-independent sequence under $\mbE$ satisfying: \\
	there exist a random variable $X$ with $C_{\V}(|X|^r)<\infty$ for some $r\in [1,2)$ and a constant $C$ such that
	\begin{equation}
		\V(|X_n|\geq x)\leq C\V(|X|\geq x), \quad \forall x\geq0, \ \forall n\geq 1.\notag
	\end{equation} 
	Then
	\begin{equation}
		(WLLN)\quad \quad  \lim_{n\rightarrow\infty}\V\left(\frac{\sum_{j=1}^{n}(X_j-\mbe[X_j])}{n^{\frac{1}{r}}}\leq -\epsilon \ or\ \frac{\sum_{j=1}^{n}(X_j-\mbE[X_j])}{n^{\frac{1}{r}}}\geq  \epsilon \right)=0,\ \forall \epsilon>0.\label{eq7}\\
	\end{equation}
\end{theorem}		

\begin{theorem}\label{th2}
	Under the same conditions of Theorem \ref{th1}, we have		
     \begin{equation}
		(SLLN)\quad \quad  \V\left(\liminf_{n\rightarrow\infty}\frac{\sum_{i=1}^{n}(X_i-\mbe[X_i])}{n^{\frac{1}{r}}}<0\ or\ \limsup_{n\rightarrow\infty}\frac{\sum_{i=1}^{n}(X_i-\mbE[X_i])}{n^{\frac{1}{r}}}>0\right)=0.\label{eq6}
	\end{equation}
\end{theorem}
	
\par We have known that Peng's independence implies the pseudo-independence, and by \eqref{eq3} and \eqref{eq4}, we can obtain the following corollary, which is similar to (3.9) of Theorem 3.4 in \cite{Zha23}. Actually, \eqref{eq7} also holds under the conditions of the following Corollary 3.1.
\begin{corollary}
	For an i.i.d. sequence $\{X_n, n\geq 1\}$ on $(\Omega, \mathcal{F},\mbE)$ satisfying $C_{\V}(|X_1|^r)<\infty$ for some $r \in [1,2)$, we have
	\begin{equation}
	\V\left(\liminf_{n\rightarrow\infty}\frac{\sum_{i=1}^{n}(X_i-\mbe[X_i])}{n^{\frac{1}{r}}}<0\ or\ \limsup_{n\rightarrow\infty}\frac{\sum_{i=1}^{n}(X_i-\mbE[X_i])}{n^{\frac{1}{r}}}>0\right)=0.\notag
	\end{equation}
\end{corollary}

\section{Proofs}
\par \ \ \ \ \ In this section, we turn to the proofs of Theorem \ref{th1} and Theorem \ref{th2}.
\subsection{Marcinkiewicz-type weak law of large numbers}
\begin{proof}[\bf{Proof of Theorem \ref{th1}}]
	Denote the truncated random variables $Y_j\triangleq(-j^{\frac{1}{r}})\vee X_j\wedge j^{\frac{1}{r}}$, \ $\forall j\in \mathbb{N}^*$, then 
	\begin{equation}
		\mbe[Y_j]\leq E_{P}[Y_j|\mathcal{F}_{j-1}]\leq \mbE[Y_j].\notag
	\end{equation}
    Note that
	\begin{equation}
		\frac{1}{n^{\frac{1}{r}}}\sum_{j=1}^{n}(X_j-\mbE[X_j])\leq \frac{1}{n^{\frac{1}{r}}}\sum_{j=1}^{n}|X_j-Y_j|+\frac{1}{n^{\frac{1}{r}}}\sum_{j=1}^{n}(Y_j-E_{P}[Y_j|\mathcal{F}_{j-1}])+\frac{1}{n^{\frac{1}{r}}}\sum_{j=1}^{n}|\mbE[Y_j]-\mbE[X_j]|.\label{eq12}
	\end{equation}
	We prove \eqref{eq7} in three steps.\\
	\par \textbf{STEP 1.} We show 
	\begin{equation}
		\lim_{n\rightarrow\infty}\frac{1}{n^{\frac{1}{r}}}\sum_{j=1}^{n}|\mbE[Y_j]-\mbE[X_j]|=0.\label{eq10}
	\end{equation}
	By the sub-additivity of $\mbE$, we have
	\begin{equation}
		|\mbE[Y_j]-\mbE[X_j]|\leq \mbE[|Y_j-X_j|]=\mbE[(|X_j|-j^{\frac{1}{r}})^+].\notag
	\end{equation}
	Notice that
	\begin{align}
	(|X_j|-j^{\frac{1}{r}})^+&=\sum_{i=j}^{\infty}(|X_j|-j^{\frac{1}{r}})^+I_{[i^{\frac{1}{r}}\leq|X_j|< (i+1)^{\frac{1}{r}}]}\notag\\
	&\leq\sum_{i=j}^{\infty}((i+1)^{\frac{1}{r}}-j^{\frac{1}{r}})^+I_{[i^{\frac{1}{r}}\leq|X_j|< (i+1)^{\frac{1}{r}}]}\notag\\
	&=\sum_{i=j}^{\infty}((i+1)^{\frac{1}{r}}-j^{\frac{1}{r}})(I_{[|X_j|\geq i^{\frac{1}{r}}]}-I_{[|X_j|\geq (i+1)^{\frac{1}{r}}]})\notag\\
	&=((j+1)^{\frac{1}{r}}-j^{\frac{1}{r}})I_{[|X_j|\geq j^{\frac{1}{r}}]}+\sum_{i=j+1}^{\infty}((i+1)^{\frac{1}{r}}-i^{\frac{1}{r}})I_{[|X_j|\geq i^{\frac{1}{r}}]}.\notag
	\end{align}
	In the case $r=1$, for $j\geq1$,
	\begin{align}
		\mbE[(|X_j|-j)^+]&\leq \V(|X_j|\geq j)+\sum_{i=j+1}^{\infty}\V(|X_j|\geq i)\notag\\
		&\leq C\V(|X|\geq j)+C\sum_{i=j+1}^{\infty}\V(|X|\geq i)\notag\\
		&\leq C\V(|X|\geq j)+C\int_{j}^{+\infty}\V(|X|\geq t)dt\rightarrow0,\ j\rightarrow\infty.\notag
	\end{align}
	Therefore, 
	\begin{equation}
		\frac{1}{n}\sum_{j=1}^{n}|\mbE[Y_j]-\mbE[X_j]|\leq \frac{1}{n}\sum_{j=1}^{n}\mbE[(|X_j|-j)^+]\rightarrow0,n\rightarrow\infty,\label{eq11}
	\end{equation}
 then \eqref{eq10} is proved in the case $r=1$.\\
	In the case $r>1$, for $j\geq1$,
	\begin{align}
		(|X_j|-j^{\frac{1}{r}})^+&\leq\frac{1}{r}\cdot j^{\frac{1}{r}-1}I_{[|X_j|\geq j^{\frac{1}{r}}]}+\sum_{i=j+1}^{\infty}\frac{1}{r}\cdot i^{\frac{1}{r}-1}I_{[|X_j|\geq i^{\frac{1}{r}]}}\notag\\
		&\leq \frac{1}{r}\cdot I_{[|X_j|\geq j^{\frac{1}{r}}]}+\frac{1}{r}\cdot\sum_{i=j+1}^{\infty} i^{\frac{1}{r}-1}I_{[|X_j|^r\geq i]}\notag.
	\end{align}
	For each $P\in \mathcal{P}$,
	\begin{align}
		E_{P}[(|X_j|-j^{\frac{1}{r}})^+]&\leq \frac{1}{r}P(|X_j|\geq j^{\frac{1}{r}})+\frac{1}{r}\sum_{i=j+1}^{\infty}i^{\frac{1}{r}-1}P(|X_j|^r\geq i)\notag\\
		&\leq \frac{1}{r}\V(|X_j|\geq j^{\frac{1}{r}})+\frac{1}{r}\sum_{i=j+1}^{\infty}i^{\frac{1}{r}-1}\V(|X_j|^r\geq i).\notag
	\end{align}
	Then taking the supremum over $P\in\mathcal{P}$, we have
	\begin{align}
		\mbE[(|X_j|-j^{\frac{1}{r}})^+]&\leq \frac{1}{r}\V(|X_j|\geq j^{\frac{1}{r}})+\frac{1}{r}\sum_{i=j+1}^{\infty}i^{\frac{1}{r}-1}\V(|X_j|^r\geq i)\notag\\
		&\leq \frac{C}{r}\V(|X|\geq j^{\frac{1}{r}})+\frac{C}{r}\sum_{i=j+1}^{\infty}i^{\frac{1}{r}-1}\V(|X|^r\geq i).\notag
	\end{align}
	It follows that
	\begin{align}
		\sum_{j=1}^{\infty}\frac{\mbE[(|X_j|-j^{\frac{1}{r}})^+]}{j^{\frac{1}{r}}}
		&\leq \frac{C}{r}\sum_{j=1}^{\infty}\V(|X|\geq j^{\frac{1}{r}})+\frac{C}{r}\sum_{j=1}^{\infty}\sum_{i=j+1}^{\infty}\frac{i^{\frac{1}{r}-1}\V(|X|^r\geq i)}{j^{\frac{1}{r}}}\notag\\
		&\leq \frac{C}{r}\cdot C_{\V}(|X|^r)+\frac{C}{r}\sum_{i=2}^{\infty}i^{\frac{1}{r}-1}\V(|X|^r\geq i)\sum_{j=1}^{i-1}j^{-\frac{1}{r}}\notag\\
		&\leq \frac{C}{r}\cdot C_{\V}(|X|^r)+\frac{C}{r}\sum_{i=2}^{\infty}i^{\frac{1}{r}-1}\V(|X|^r\geq i)\cdot \frac{r}{r-1}i^{1-\frac{1}{r}}\notag\\
		&\leq \frac{2C}{r-1}C_{\V}(|X|^r)<\infty,\notag
	\end{align} 
	where the third inequality is due to the fact that for $i\geq 2$,
	\begin{equation}
		\sum_{j=1}^{i-1}j^{-\frac{1}{r}}\leq \int_{0}^{i}x^{-\frac{1}{r}}dx=\frac{r}{r-1}i^{1-\frac{1}{r}}.\notag
	\end{equation}
	By Kronecker's lemma, we have
	\begin{equation}
		\frac{\sum_{j=1}^{n}\mbE[(|X_j|-j^{\frac{1}{r}})^+]}{n^{\frac{1}{r}}}\rightarrow0,n\rightarrow\infty.\notag
	\end{equation}  
	\eqref{eq10} is proved.
	 \par \textbf{STEP 2.} For $\epsilon>0$, we show
	 \begin{equation}
	 	P\left(\frac{1}{n^{\frac{1}{r}}}\sum_{j=1}^{n}(Y_j-E_{P}[Y_j|\mathcal{F}_{j-1}])\geq\frac{\epsilon}{4}\right)\rightarrow0,n\rightarrow\infty, \text{uniformly for all}\ P\in \mathcal{P}.\notag
	 \end{equation} 
	 By Chebyshev's inequality, for all $\epsilon>0$ and $P\in \mathcal{P}$ we have
	 \begin{equation}
	 	P\left(\frac{1}{{n^{\frac{1}{r}}}}\sum_{j=1}^{n}(Y_j-E_{P}[Y_j|\mathcal{F}_{j-1}])\geq\frac{\epsilon}{4}\right)\leq\frac{16}{\epsilon^2}\cdot\frac{\sum_{j=1}^{n}E_{P}[Y_j^2]}{n^{\frac{2}{r}}}\leq\frac{16}{\epsilon^2}\cdot\frac{\sum_{j=1}^{n}\mbE[Y_j^2]}{n^{\frac{2}{r}}}. \notag
	 \end{equation}
	 By Kronecker's lemma, we only need to prove
	 \begin{equation}
	 	\sum_{j=1}^{\infty}\frac{\mbE[Y_j^2]}{j^{\frac{2}{r}}}<\infty ,\notag
	 \end{equation}
	 which yields
	 \begin{equation}
	 	\lim_{n\rightarrow\infty}\frac{\sum_{j=1}^{n}\mbE[Y_j^2]}{n^{\frac{2}{r}}}=0.\notag
	 \end{equation}
	  Note that
	 \begin{align}
	 	Y_j^2
	 	&=\sum_{i=1}^{j}X_j^2I_{[(i-1)^{\frac{1}{r}}\leq |X_j|<i^{\frac{1}{r}}]}+j^{\frac{2}{r}}I_{[|X_j|\geq j^{\frac{1}{r}}]}\notag\\
	 	&\leq \sum_{i=1}^{j}i^{\frac{2}{r}}(I_{[|X_j|\geq (i-1)^{\frac{1}{r}}]}-I_{[|X_j|\geq i^{\frac{1}{r}}]})+j^{\frac{2}{r}}I_{[|X_j|\geq j^{\frac{1}{r}}]}\notag\\
	 	&\leq 1+\sum_{i=1}^{j-1}\frac{2}{r}(i+1)^{\frac{2}{r}-1}I_{[|X_j|\geq i^{\frac{1}{r}}]}\notag\\
	 	&\leq 1+4\sum_{i=1}^{j}i^{\frac{2}{r}-1}I_{[|X_j|\geq i^{\frac{1}{r}}]},\notag
	 \end{align}
	 then 
	 \begin{align}
	 	\mbE[Y_j^2]&\leq 1+4\sum_{i=1}^{j}i^{\frac{2}{r}-1}\V(|X_j|\geq i^{\frac{1}{r}})\notag\\
	 	&\leq 1+4C\sum_{i=1}^{j}i^{\frac{2}{r}-1}\V(|X|\geq i^{\frac{1}{r}}).\notag
	 \end{align}
	 It follows that 
	 \begin{align}
	 	\sum_{j=1}^{\infty}\frac{\mbE[Y_j^2]}{j^{\frac{2}{r}}}&\leq \sum_{j=1}^{\infty}\frac{1}{j^{\frac{2}{r}}}+4C\sum_{j=1}^{\infty}\sum_{i=1}^{j}\frac{i^{\frac{2}{r}-1}\V(|X|\geq i^{\frac{1}{r}})}{j^{\frac{2}{r}}}\notag\\
	 	&=\sum_{j=1}^{\infty}\frac{1}{j^{\frac{2}{r}}}+4C\sum_{i=1}^{\infty}i^{\frac{2}{r}-1}\V(|X|\geq i^{\frac{1}{r}})\sum_{j=i}^{\infty}\frac{1}{j^{\frac{2}{r}}}\notag\\
	 	&\leq (1+4C)\sum_{j=1}^{\infty}\frac{1}{j^{\frac{2}{r}}}+4C\sum_{i=2}^{\infty}i^{\frac{2}{r}-1}\V(|X|\geq i^{\frac{1}{r}})\sum_{j=i}^{\infty}\frac{1}{j^{\frac{2}{r}}}\notag\\
	 	&\leq (1+4C)\sum_{j=1}^{\infty}\frac{1}{j^{\frac{2}{r}}}+4C\sum_{i=2}^{\infty}i^{\frac{2}{r}-1}\V(|X|\geq i^{\frac{1}{r}})\cdot\frac{r}{2-r}(i-1)^{-\frac{2}{r}+1}\notag\\
	 	&\leq (1+4C)\sum_{j=1}^{\infty}\frac{1}{j^{\frac{2}{r}}}+\frac{8rC}{2-r}C_{\V}(|X|^r)<\infty,\notag
	 \end{align}
	 where we've used the fact that for $i\geq 2$,
	 \begin{equation}
	 	\sum_{j=i}^{\infty}j^{-\frac{2}{r}}\leq \int_{i-1}^{\infty}x^{-\frac{2}{r}}dx=\frac{r}{2-r}(i-1)^{-\frac{2}{r}+1}.\notag
	 \end{equation}
	 \par \textbf{STEP 3.} For $\epsilon>0$, we show
	 \begin{equation}
	 	P\left(\frac{1}{n^{\frac{1}{r}}}\sum_{j=1}^{n}|X_j-Y_j|\geq\frac{\epsilon}{4}\right)\rightarrow0,n\rightarrow\infty, \text{uniformly for all}\ P\in \mathcal{P}.\notag
	 \end{equation}
	  Similar to the second step, we use Markov's inequality to get that for all $\epsilon>0$ and $P\in \mathcal{P}$,
	  \begin{equation}
	 	P\left(\frac{1}{{n^{\frac{1}{r}}}}\sum_{j=1}^{n}|X_j-Y_j|\geq\frac{\epsilon}{4}\right)\leq\frac{4}{\epsilon}\cdot\frac{\sum_{j=1}^{n}E_{P}[|X_j-Y_j|]}{n^{\frac{1}{r}}}\leq\frac{4}{\epsilon}\cdot\frac{\sum_{j=1}^{n}\mbE[|X_j-Y_j|]}{n^{\frac{1}{r}}}.\notag
	 \end{equation}
	 It has been proved in the first step that
	 \begin{equation}
	 	\sum_{j=1}^{\infty}\frac{\mbE[|X_j-Y_j|]}{j^{\frac{1}{r}}}=\sum_{j=1}^{\infty}\frac{\mbE[(|X_j|-j^{\frac{1}{r}})^+]}{j^{\frac{1}{r}}}<\infty.\notag
	 \end{equation}
	 Then by Kronecker's lemma, for all $P\in \mathcal{P}$,
	 \begin{equation}
	 	P\left(\frac{1}{{n^{\frac{1}{r}}}}\sum_{j=1}^{n}|X_j-Y_j|\geq\frac{\epsilon}{4}\right)\leq \frac{4}{\epsilon}\cdot\frac{\sum_{j=1}^{n}\mbE[|X_j-Y_j|]}{n^{\frac{1}{r}}}\rightarrow0,\ n\rightarrow\infty.\notag
	 \end{equation}
	 By \eqref{eq12}, we obtain that for all $\epsilon>0$ and $P \in \mathcal{P},$
	 \begin{align}
	 	P\left(\frac{1}{{n^{\frac{1}{r}}}}\sum_{j=1}^{n}(X_j-\mbE[X_j])\geq \epsilon\right)&\leq P\left(\frac{1}{{n^{\frac{1}{r}}}}\sum_{j=1}^{n}|X_j-Y_j|\geq\frac{\epsilon}{4}\right)+ P\left(\frac{1}{{n^{\frac{1}{r}}}}\sum_{j=1}^{n}(Y_j-E_{P}[Y_j|\mathcal{F}_{j-1}])\geq\frac{\epsilon}{4}\right)\notag\\
	 	&+P\left(\frac{1}{n^{\frac{1}{r}}}\sum_{j=1}^{n}|\mbE[Y_j]-\mbE[X_j]|\geq \frac{\epsilon}{2}\right)\notag\\
	 	&\leq \frac{4}{\epsilon}\cdot\frac{\sum_{j=1}^{n}\mbE[|X_j-Y_j|]}{n^{\frac{1}{r}}}+ \frac{16}{\epsilon^2}\cdot\frac{\sum_{j=1}^{n}\mbE[Y_j^2]}{n^{\frac{2}{r}}}\notag\\
	 	&+\frac{2}{\epsilon}\cdot\frac{1}{n^{\frac{1}{r}}}\sum_{j=1}^{n}|\mbE[Y_j]-\mbE[X_j]|\rightarrow0,n\rightarrow\infty.\notag
	 \end{align}
	 Then taking the supremum over $P\in\mathcal{P}$, we have
	 \begin{equation}
	 	\lim_{n\rightarrow\infty}\V\left(\frac{\sum_{j=1}^{n}(X_j-\mbE[X_j])}{n^{\frac{1}{r}}}\geq  \epsilon \right)=0,\ \forall \epsilon>0.\notag
	 \end{equation}
	 On the other hand, 
	 \begin{equation}
	 	\lim_{n\rightarrow\infty}\V\left(\frac{\sum_{j=1}^{n}(X_j-\mbe[X_j])}{n^{\frac{1}{r}}}\leq-\epsilon \right)=\lim_{n\rightarrow\infty}\V\left(\frac{\sum_{j=1}^{n}(-X_j-\mbE[-X_j])}{n^{\frac{1}{r}}}\geq  \epsilon \right)=0,\ \forall \epsilon>0.\notag
	 \end{equation}
	 By the sub-additivity of $\V$, \eqref{eq7} is proved. The proof is now completed.
\end{proof}
\subsection{Marcinkiewicz-type strong law of large numbers}
	 \par Now, we show the proof of \eqref{eq6}.
\begin{proof}[\bf{Proof of Theorem \ref{th2}}]
Denote the truncated random variables $Y_j\triangleq(-j^{\frac{1}{r}})\vee X_j\wedge j^{\frac{1}{r}}$, \ $\forall j\in \mathbb{N}^*$, which is the same as in the proof of Theorem \ref{th1}.\\
Note that 
\begin{equation}
	\frac{1}{n^{\frac{1}{r}}}\sum_{j=1}^{n}(X_j-\mbE[X_j])\leq \frac{1}{n^{\frac{1}{r}}}\sum_{j=1}^{n}|X_j-Y_j|+\frac{1}{n^{\frac{1}{r}}}\sum_{j=1}^{n}(Y_j-\mbE[Y_j])+\frac{1}{n^{\frac{1}{r}}}\sum_{j=1}^{n}|\mbE[Y_j]-\mbE[X_j]|.\label{eq13}
\end{equation}
Similar to the proof of Theorem \ref{th1}, we prove \eqref{eq6} in three steps as well.
\par \textbf{STEP 1.} We show 
\begin{equation}
	\V\left(\limsup_{n\rightarrow\infty}\frac{1}{n^{\frac{1}{r}}}\sum_{j=1}^{n}(Y_j-\mbE[Y_j])>0\right)=0.\label{eq14}
\end{equation}
Note that $\{\sum_{j=1}^{n}\frac{Y_j-E_{P}[Y_j|\mathcal{F}_{j-1}]}{j^{\frac{1}{r}}},\mathcal{F}_{n},n\geq 1\}$ is a martingale for each $P\in \mathcal{P}$. By a simple calculation,
\begin{equation}
	E_{P}\left[\sum_{j=1}^{\infty}E_P\left[\frac{(Y_j-E_{P}[Y_j|\mathcal{F}_{j-1}])^2}{j^{\frac{2}{r}}}\Bigg|\mathcal{F}_{j-1}\right]\right]\leq\sum_{j=1}^{\infty}\frac{E_P[Y_j^2]}{j^{\frac{2}{r}}}\leq \sum_{j=1}^{\infty}\frac{\mbE[Y_j^2]}{j^{\frac{2}{r}}}<\infty,\notag
\end{equation}  
where the last inequality has been proved in the second step of the proof of Theorem \ref{th1}.\\
Thus, 
\begin{equation}
	P\left(\sum_{j=1}^{\infty}E_P\left[\frac{(Y_j-E_{P}[Y_j|\mathcal{F}_{j-1}])^2}{j^{\frac{2}{r}}}\Bigg|\mathcal{F}_{j-1}\right]<\infty\right)=1\notag
\end{equation}
implies
\begin{equation}
	P\left(\sum_{j=1}^{\infty}\frac{Y_j-E_{P}[Y_j|\mathcal{F}_{j-1}]}{j^{\frac{1}{r}}}<\infty\right)=1\notag
\end{equation}
by Lemma \ref{le1}.
\par By Kronecker's lemma, we get
\begin{equation}
	\lim_{n\rightarrow\infty}\frac{1}{n^{\frac{1}{r}}}\sum_{j=1}^{n}(Y_j-E_P[Y_j|\mathcal{F}_{j-1}])=0,\ P-a.s.,\ \forall P\in \mathcal{P}.\notag
\end{equation}
Notice that
\begin{equation}
	P\left(\limsup_{n\rightarrow\infty}\frac{1}{n^{\frac{1}{r}}}\sum_{j=1}^{n}(Y_j-\mbE[Y_j])>0\right)\leq P\left(\limsup_{n\rightarrow\infty}\frac{1}{n^{\frac{1}{r}}}\sum_{j=1}^{n}(Y_j-E_P[Y_j|\mathcal{F}_{j-1}])>0\right)=0, \forall P\in \mathcal{P}.\notag
\end{equation}
Taking the supremum over $P\in \mathcal{P}$, we get \eqref{eq14}.
\par \textbf{STEP 2.} We show 
\begin{equation}
	\V\left(\limsup_{n\rightarrow\infty}\frac{1}{n^{\frac{1}{r}}}\sum_{j=1}^{n}|X_j-Y_j|>0\right)=0.\label{eq15}
\end{equation}
It's clear that 
\begin{equation}
	\sum_{j=1}^{\infty}\V(|X_j-Y_j|>0)=\sum_{j=1}^{\infty}\V(|X_j|\geq j^{\frac{1}{r}})\leq C\sum_{j=1}^{\infty}\V(|X|\geq j^{\frac{1}{r}})\leq C\cdot C_{\V}(|X|^r)< \infty.\notag
\end{equation}
By the Borel-Cantelli lemma, we have
\begin{equation}
	\V\left(\limsup_{n\rightarrow\infty}\frac{1}{n^{\frac{1}{r}}}\sum_{j=1}^{n}|X_j-Y_j|>0\right)\leq \V(|X_j-Y_j|>0,i.o.)=0.\notag
\end{equation}
\par \textbf{STEP 3.} It's true that
\begin{equation}
	\lim_{n\rightarrow\infty}\frac{1}{n^{\frac{1}{r}}}\sum_{j=1}^{n}|\mbE[Y_j]-\mbE[X_j]|=0,\label{eq16}
\end{equation}
which has been proved in the first step of the proof of Theorem \ref{th1}.\\
By \eqref{eq14},\eqref{eq15} and \eqref{eq16}, we have
\begin{align}
	\V\left(\limsup_{n\rightarrow\infty}\frac{1}{n^{\frac{1}{r}}}\sum_{j=1}^{n}(X_j-\mbE[X_j])>0\right)&\leq \V\left(\limsup_{n\rightarrow\infty}\frac{1}{n^{\frac{1}{r}}}\sum_{j=1}^{n}|X_j-Y_j|>0\right)\notag\\
	&+\V\left(\limsup_{n\rightarrow\infty}\frac{1}{n^{\frac{1}{r}}}\sum_{j=1}^{n}(Y_j-\mbE[Y_j])>0\right)=0.\notag
\end{align}
On the other hand,
\begin{equation}
	\V\left(\liminf_{n\rightarrow\infty}\frac{\sum_{i=1}^{n}(X_i-\mbe[X_i])}{n^{\frac{1}{r}}}<0\right)=\V\left(\limsup_{n\rightarrow\infty}\frac{\sum_{i=1}^{n}(-X_i-\mbE[-X_i])}{n^{\frac{1}{r}}}>0\right)=0.\notag
\end{equation}
The proof is completed.
\end{proof}

\bibliographystyle{plain}

\end{document}